\documentclass[a4paper,12pt,leqno]{amsart}
\usepackage[latin1]{inputenc}
\usepackage[T1]{fontenc}
\usepackage{amssymb,calrsfs,ae}
\usepackage{hyperref}

\swapnumbers
\theoremstyle{plain}
\newtheorem{prop}{Proposition}[section]
\newtheorem{theo}[prop]{Theorem}
\newtheorem{cor}[prop]{Corollary}

\theoremstyle{definition}
\newtheorem{defi}[prop]{Definition}

\theoremstyle{remark}

\newcommand{\thismonth}{\ifcase\month\or
  January\or February\or March\or April\or May\or June\or July\or
  August\or September\or October\or November\or December\fi
  \space\number\year}
\DeclareMathAlphabet{\mathrmsl}{OT1}{cmr}{m}{sl}

\newcommand{\oper}[3][n]{\newcommand{#2}{\mathop{\mathrm{#3}}%
\ifx n#1\nolimits\else\limits\fi} }
\newcommand{\rsoper}[3][n]{\newcommand{#2}{\mathop{\mathrmsl{#3}}%
\ifx n#1\nolimits\else\limits\fi} }
\newcommand{\NM}{{\mathbb N}}

\newcommand{\ZM}{{\mathbb Z}}

\newcommand{\vol}{\operatorname{vol}}
\newcommand{\Scal}{\operatorname{Scal}}

\newcommand{\csum}{\mathfrak{S}}

\newcommand{\tv}{\mathcal{T}}

\newcommand{\II}{\mathbb{I}}

\newcommand{\cerc}{{\mathbb S}}

\renewcommand{\exp}{\operatorname{e}}
\newcommand{\gb}{\overline{g}}
\newcommand{\bg}{\gb}
\newcommand{\Id}{\operatorname{Id}}

\newcommand{\proofof}[1]{\end{#1}\begin{proof}}

\newcounter{mnotecount}[section]
\renewcommand{\themnotecount}{\thesection.\arabic{mnotecount}}
\newcommand{\mnote}[1]
{\protect{\stepcounter{mnotecount}}$^{\mbox{\footnotesize  $
      \bullet$\themnotecount}}$ \marginpar{\raggedright\tiny\em
    $\!\!\!\!\!\!\,\bullet$\themnotecount: #1} }

\begin{document}

\title[CR Seifert 3-manifolds]{An invariant of Cauchy-Riemann Seifert 
3-manifolds and applications}
\author{Olivier Biquard}
\address{Institut de Recherche Mathématique Avancée\\
UMR 7501 CNRS -- Université Louis Pasteur \\ Strasbourg (France)}
\email{olivier.biquard@math.u-strasbg.fr}
\author{Marc Herzlich}
\address{Institut de Mathématiques et Modélisation de Montpellier\\ 
UMR 5149 CNRS -- Université Mont­pellier~II\\ France}
\email{herzlich@math.univ-montp2.fr}

\begin{abstract}
We compute a recently introduced geometric invariant of stricly
pseudoconvex CR 3-manifolds for certain circle invariant
spherical CR structures on Seifert manifolds.
We give applications to the problem of filling the CR manifold by a
complex hyperbolic manifold, and more generally by a Kähler-Einstein
or an Einstein metric.
\end{abstract}
\thanks{The second author is supported in part by the Young Researchers 
{\sc aci} program of the French Ministry of Research.}
\maketitle
\section{Introduction.}

In \cite{ob-mh1} we introduced a new invariant, called the $\nu$-invariant, of
strictly pseudoconvex Cauchy-Riemann (CR) compact 3-manifolds.
This invariant is an analogue in CR geometry of the $\eta$-invariant in
conformal geometry. The definition of the $\nu$-invariant is rather abstract 
and makes it difficult to get explicit expressions. The aim of this paper is 
to provide a computation of $\nu$ in the simple, yet interesting, case
of certain spherical $\cerc^1$-invariant CR structures on Seifert 
manifolds, and to deduce some geometric applications.

\medskip
The spherical (that is, locally isomorphic to the standard CR
3-sphere) CR structures we are interested in come with an action of
$\cerc^1$ without fixed point and transverse to the contact distribution
(see \cite{ka-tsu} for a classification of all spherical
$\cerc^1$-invariant CR structures). They appear as
orbifold $\cerc^1$-bundles over $2$-dimensional orbifolds. 
At each orbifold point, the orbifold data consists of the following:
the local fundamental group is $\ZM/\alpha\ZM$ ($\alpha\in\NM^*$), and a generator
acts on a local chart around $p$ of the basis manifold as 
$\exp^{i\frac{2\pi\beta}{\alpha}}$ and on the fiber as 
$\exp^{i\frac{2\pi\gamma}{\alpha}}$ with $\beta$ and $\gamma$ prime to $\alpha$.
The orbifold $\cerc^1$ bundle is topologically classified by the degree (first Chern
number) $d$ of the bundle---in this case a rational number.
We then endow the manifold with an invariant strictly pseudoconvex CR structure: 
The underlying contact structure is provided by a constant curvature equivariant
connection $1$-form on the bundle, whereas the complex structure is induced from 
the basis Riemann surface. The pseudoconvexity condition constrains the degree $d$ 
to be negative. Our main result then reads:

\begin{theo}\label{maintheo}
Let $X$ be a compact spherical Cauchy-Riemann $3$-manifold which is 
a $\cerc^1$ orbifold bundle of degree $d<0$ over a compact orbifold Riemann
surface $\Sigma$ of Euler characteristic $\chi$ \upn{(}a rational number\upn{)}.
The $\nu$-invariant of $X$ is 
\begin{equation}\label{eq:nu} 
\nu(X) = - d - 3 - \frac{\chi^2}{4d} - 12 \sum_{j=1}^{p} 
s(\alpha_j,\beta_j,\gamma_j)
\end{equation}
where $s(\alpha,\beta,\gamma)$ is the Dedekind sum
$\frac{1}{4\alpha} \sum\limits_{k=1}^{\alpha -1}
\cot \left( \frac{k\beta\pi}{\alpha}\right)
\cot \left( \frac{k\gamma\pi}{\alpha}\right)$ .
\end{theo}

\medskip

We first remark that for this class of CR manifolds,
the $\nu$-invariant depends only on the topology, and not, for instance,
on the basis complex structure. This is {\it a priori}
known, since the gradient of $\nu$ is the Cartan curvature 
\cite[Theorem 8.1]{ob-mh1}, which vanishes for spherical CR manifolds. 

\medskip

We now pass to geometric applications. Such $X$ being
locally isomorphic to the standard CR sphere $\cerc^3$,
it is the boundary at infinity of a complex hyperbolic
metric defined in a neighbourhood $(0,\varepsilon]\times X$ of $X$
(in the case of the $3$-sphere we get the Bergmann metric on the 4-ball).

From \cite[Theorem 1.2]{ob-mh1}, we get the following obstruction for 
this neighbourhood to have a global extension to a smooth complex 
hyperbolic surface (with only one end):

\begin{cor}\label{cor:1}
  Let $X$ be as in Theorem \ref{maintheo}. If $X^3$ is the boundary at
  infinity of a complex hyperbolic metric defined on the interior $M$ 
  of a smooth compact manifold ${\bar M}^4$ with boundary $X$,
  then one has necessarily $ \nu(X) = - \chi(\bar{M}) + 3 \tau(\bar{M})$,
  where $\chi(\bar{M})$ and $\tau(\bar{M})$ denote the Euler characteristic 
  and signature of $\bar{M}$. In particular, $\nu(X)$, as given by
  the formula \upn{(}\ref{eq:nu}\upn{)}, is an integer.
\end{cor}

We can restate this in a special case :

\begin{cor}\label{cor:2}
  Let $X$ be a $\cerc^1$-bundle of degree $d$ over a Riemann surface $\Sigma$ 
  of Euler characteristic $\chi$, with a $\cerc^1$-invariant 
  spherical CR structure.
  If $ \frac{ \chi^2}{4d}$ is not an integer then $X$ is not
  the boundary at infinity of a complex hyperbolic metric.
\end{cor}

The case $d=\frac{\chi}{2}$ yields an integer, and indeed, if $\Sigma$
is hyperbolic, $\bar{M}$ can be taken to be the disk bundle of a 
square root of the tangent bundle of $\Sigma$, which is well known to carry 
a complex hyperbolic metric issued from a representation of $\pi_1(\Sigma)$ 
in $SU(1,1) \subset SU(1,2)$. Our obstruction then gives an interesting 
hint on whether a CR flat Seifert $3$-manifold may appear as a quotient of the 
complement of the limit set in the $3$-sphere of some discrete fixed 
point-free subgroup of $SU(1,2)$ \cite{apanasov}.

More generally, the calculation in Theorem \ref{maintheo} gives an
obstruction for $X$ to be the boundary at infinity of a
Kähler-Einstein metric. The manifolds considered in this paper are
known to bound a complex Stein space with at most a finite number of
singular points \cite{harvey-lawson} and one may wish to endow it
with a K\"ahler-Einstein metric as in Cheng-Yau \cite{cheng-yau}. 
The type of metric to be considered has
the same kind of asymptotic expansion near the boundary $X$ as the
Bergmann metric \cite{biquard-symmetric}; we called them ``asymptotically 
complex hyperbolic'' (ACH) in \cite{ob-mh1}. If no singular points are
present and if the Cheng-Yau metric exists, one gets 
from the Miyaoka-Yau inequality the following:

\begin{cor}\label{cor:3}
  Let $X$ be as in Theorem \ref{maintheo}. If $X$ is the boundary at
  infinity of an ACH Kähler-Einstein metric on $M^4$, then
  $$ \chi(\bar{M}) - 3 \tau(\bar{M}) \geqslant -\nu(X) =
  d + 3 + \frac{\chi^2}{4d} + 12 \sum_{j=1}^{p} s(\alpha_j,\beta_j,\gamma_j).$$
\end{cor}

This is a topological constraint on a filling. For more information on Stein 
fillings, see \cite{lisca-matic,stipsicz1}. We mention that this corollary 
remains true more generally for ACH Einstein (not K\"ahler) fillings, under 
the additional assumption that a Kronheimer-Mrowka contact invariant 
\cite{kro-mro} of $(\bar{M},X)$ is nonzero.
This condition is fulfilled if $M$ admits a symplectic form compatible with
the contact structure on its boundary. 

\medskip

From \cite[Theorem 5.12]{carron-pedon}, one knows that pseudoconvex
complex hyperbolic surfaces $\bar{M}$ have vanishing third homology group 
$H_3(\bar{M},\mathbb{Z})$. Hence no multiple ends can occur, but 
one expects orbifold singularities or cusps to appear in the interior
of a complex hyperbolic 
filling. The complex hyperbolic cusps can be compactified to yield a complex 
orbifold surface that we note again $\bar M$,
by adding at the infinity of each cusp a quotient $\Sigma_i$ of a 2-torus.
The corollaries \ref{cor:1} and \ref{cor:3}
remain true in this case, with the Euler characteristic and the
signature of $\bar M$ replaced by their orbifold versions:
In case $\ell$ cusps are present, there is an additional contribution in the
signature coming from the self-intersection of each 2-torus at infinity. 
Namely, one has to
consider the modified signature \cite[proposition 3.4]{biquard-rigidite-fini}
$$ \tau_{\textrm cusp}(\bar{M})=\tau(\bar{M}) - \frac{1}{3}\sum_1^\ell
[\Sigma_i]\cdot[\Sigma_i] . $$
Of course, corollary \ref{cor:2}  is no more true, since the
characteristic numbers are now rational; the denominator of $\nu$ only
gives an hint on the order of the singularities needed to fill $X$.

Another important point is to compare these results with those obtained by use
of the Burns-Epstein $\mu$-invariant \cite{BE88,BE90b} (it is already suggested 
at the end of \cite{BE90b} that obstructions follow from computations of $\mu$).
The $\mu$ invariant is defined on strictly pseudoconvex CR 3-manifolds with
trivial tangent holomorphic bundle only. Roughly speaking, it comes from
Chern-Simons-type constructions (integration of a local formula), whereas the 
$\nu$-invariant is extracted from the Atiyah-Patodi-Singer $\eta$-invariant.
The relation between $\mu$ and $\nu$ is similar to that between the
$\eta$ and the Chern-Simons invariants:
More precisely \cite[Theorem 1.3]{ob-mh1}, when $\mu$ is defined, then
for a CR structure $J$ one has
$$ \nu(J) \ = \ 3\,\mu(J) \ + \ \text{constant} , $$
with the constant depending only of the underlying contact structure.
Burns-Epstein's version of Miyaoka-Yau \cite{BE90b} then reads, if
$X$ is the boundary at infinity of a Kähler-Einstein $M$:
\begin{equation}\label{eq:BE}
\chi(\bar{M})-\frac{1}{3}\,\bar{c}_1(\bar{M})^2 \geqslant - \mu(X) ,
\end{equation}
with equality if the metric is complex hyperbolic;
here $\bar{c}_1$ is a lift in $H^2(M,X)$ of $c_1(M)$.

A first important difference here is that our obstruction in 
Corollary \ref{cor:3} (filling by an ACH K\"ahler-Einstein metric) is
purely topological, whereas (\ref{eq:BE}) involves the complex structure.

Another important fact to be noticed, at least in the case when there is no
orbifold quotient, is that the obstructions obtained by both methods are
different: if $X$ is a 
$\cerc^1$-bundle over the Riemann surface $\Sigma$, then the $\mu$-invariant, 
being defined by a local formula, is multiplicative on finite coverings 
\cite{BE88,BE90b}. Hence the values are
\begin{equation}\label{eq:mu.vs.nu} 
\mu = \frac{\chi^2}{4d} 
\ \textrm{ whereas } \  \nu = - \frac{ \chi^2}{4d} - d - 3.
\end{equation}
Equation (\ref{eq:BE}) implies that $3\mu$ must be an integer, {\it i.e.}
$\frac{3\chi^2}{4d}$ must belong to $\mathbb{Z}$, a condition that is 
weaker than corollary \ref{cor:2}, by a factor $3$.

\medskip

The reader will find the computation of $\nu$ in the next section. The beginning
of section $3$ is then devoted to a full explicitation of these results in an 
interesting particular case: that of lens spaces; formulas are given in 
Proposition 3.1. The paper then ends with the proof of the corollaries.

\bigskip

\section{Computation of the invariant}

We first describe our $3$-dimensional compact strictly pseudoconvex CR manifold 
$X$ in greater detail. The basis compact Riemannian orbifold $\Sigma$ is endowed
with a constant curvature metric and a complex structure. For any $\cerc^1$-bundle 
$X$ over it of degree $d<0$,
 there exists an equivariant connection $1$-form $\bar\theta$
of constant curvature (relative to the basis metric). We shall take 
\[ \theta \ = \ \frac{1}{2}\,\bar\theta \]
as the contact form on $X$ (the choice of the factor $\frac{1}{2}$, although 
slightly awkward at some places below, is made to ensure full consistency 
with the conventions
of our previous work \cite{ob-mh1}); in this convention, fibers then have length 
$\pi$ (rather than the more natural $2\pi$).
The complex structure $J$ on the basis can be lifted on the contact ditribution 
$H$, and the bilinear form $\gamma = d\theta(J\cdot,\cdot)$ is an invariant  
metric on $H$ under the circle action, which projects downwards as a constant 
multiple of the basis metric. From our normalization of contact form above, we have
\[ \int_{\Sigma} d\theta = -\pi d,\ \ \int_X \theta\wedge d\theta = - \pi^2 d. \]
The volume of the basis $\Sigma$ endowed with the projected metric $\gamma$ is 
$V=-\pi d$ and its curvature is $R= - \frac{2\chi}{d}$.

\medskip

One then defines a sequence of Riemannian metrics
\[ g_{\rho} = 4\rho^2\,\theta^2 + \gamma \]
on $X$, and one 
chooses an asymptotically complex hyperbolic K\"ahler-Einstein 
metric $\bg$ (with scalar curvature $-6$) on $M=[r_0,+\infty[\times X$,
defined as follows: the complex structure is trivially extended as
\[ J_{|H} = J \ \textrm{ on } H, \ \ J\partial_r = \exp^{-r}\xi ,\]
if $\xi $ is the (vertical) Reeb field of $\theta$. The K\"ahler form is 
(see \cite{ob-mh1}):
\begin{equation*}
\bar\omega = \exp^r \left( dr\land\eta + d\eta\right) - \frac{R}{2}\, d\eta
 - \frac{R^2}{12}\exp^{-r} \left(dr\land\eta - d\eta\right) + 
o(\exp^{-2r}) . \end{equation*}
It is explained in \cite{ob-mh1} why lower order terms in $\bar\omega$ are
irrelevant in all what concerns the $\nu$-invariant to be defined below.
Let us denote its $(3,1)$ curvature tensor by $R$, the
associated operator on $2$-forms by $R^{op}$ and denote by $\bar{\textrm{n}}$ and
$\II_r$ the outer unit normal and second fundamental
form of the slice $\{r\}\times X$. We also need to define for a tensor $F$ in 
$\otimes^3T^*M$, 
\[ \csum(F)(X,Y,Z) = F(X,Y,Z) + F(Y,Z,X) + F(Z,X,Y), \]
and, if $\rho$ is a $3$-form and $\sigma_1$, $\sigma_2$, $\sigma_3$
are three vectors,
\begin{equation*}
\tv\left(\rho\otimes(\sigma_1\otimes\sigma_2\otimes\sigma_3)\right) = 
d\!\vol_{\{r\}\times X} \left(\sigma_1,\sigma_2,\sigma_3\right) \,\rho .
\end{equation*} 
This definition immedia\-tely extends to a $3$-form with values 
in $\otimes^3 TM$.
 
\medskip

According to \cite{ob-mh1}, the $\nu$-invariant 
is obtained by taking the limit as $r$ goes to infinity
of the boundary contribution at $\{r\}\times X$ of the Atiyah-Patodi-Singer
expression for the characteristic number $\chi -3\tau$ of $[r_0,r]\times X$ 
with respect to the metric $\bg$. This is explicited in the 

\begin{defi} The $\nu$-invariant of $X$ is 
\[ \nu(X) \ = \ \lim_{r\rightarrow +\infty} \nu(\bg,r) \ = \ 
\lim_{r\rightarrow +\infty} \{ B(\bg,r) - 3\,\eta( 4\rho(r)^2\theta^2 + \gamma)
\}\ ,\]
where $\rho(r)= \frac{1}{4}\exp^r(1 + \frac{R}{2}\exp^{-r} + 
\frac{R^2}{12}\exp^{-2r})$, and $B(\bg,r)$ is the integral term 
\[ B(\bg,r) = -\frac{1}{12\pi^2} 
\int_{\{r\}\times X} \tv(\II_r\wedge\II_r\wedge\II_r)
  + 3\, \tv(\II_r\wedge R^{op}) - 3\, \csum (\II_r(.,R(.,.)\bar{\textrm{n}}) .
 \]
It is shown in \cite{ob-mh1} that the limit converges and gives indeed rise to a CR
invariant of $X$.
\end{defi}

\medskip

The $\eta$-invariant of (orbifold) circle bundles over (orbifold) Riemannian 
surfaces has been 
computed by Komuro \cite{komuro} and more generally by Ouyang \cite{ouyang}.
In our conventions, their results read:

\begin{theo}[Ouyang]\label{th:ouy}
The $\eta$-invariant of the metric $4\rho^2\theta^2 + \gamma$ on $X$ is equal to
\[ \frac{1}{3} \left( d + 3 + 2d \left( \frac{\pi\rho^2}{V}\chi - 
\frac{\pi^2\rho^4}{V^2}d^2  \right)\right) + 4 \sum_{j=1}^{p} 
s(\alpha_j,\beta_j,\gamma_j) . \]
\end{theo}

In our setting, this yields:
\[ 3\,\eta(r) = - \frac{d}{8}\exp^{2r} - \frac{\chi}{4}\exp^r + 
\frac{7}{24}\,\frac{\chi^2}{d} + d + 3 + 12 \sum_{j=1}^{p} 
s(\alpha_j,\beta_j,\gamma_j) . \]
The local integral terms can now be computed by using the same techniques as
in \cite{ob-mh1} and in \cite{mh-vol}. These are lengthy but otherwise
straightforward computations, and we give a few intermediate steps below.

The second fundamental forms are, up to order $2$ terms,
\[ \II_r = \left( 1+ \frac{R^2}{8}\exp^{-2r}\right) \Id_{\xi}
+ \frac{1}{2}\left( 1 + \frac{R}{2}\exp^{-r} + \frac{R^2}{8}\exp^{-2r}\right)
\Id_H + o(\exp^{-2r}) \]
(with $\xi$ the Reeb form of the contact distribution)
and the terms involving only $\II_r$ in the definition of $B(\bg,r)$ 
are easily obtained from it; one gets the
contribution:
\[ \left( \frac{3}{2}\exp^{2r} + \frac{3}{4}R \exp^r  + \frac{R^2}{4}\right)
\theta\wedge d\theta .\]
The curvature of the K\"ahler metric is formed from a constant term whose
coefficient are exactly those of complex hyperbolic curvature in a $\bg$-orthonormal
basis on $[r_0,+\infty[\times X$ such as 
\begin{multline*}
\big(\partial_r, \exp^{-r}(1 -\frac{R^2}{12}\exp^{-2r})^{-\frac{1}{2}}\xi,\\
\exp^{-r/2}(1-\frac{R}{2}\exp^{-r} + \frac{R^2}{12}\exp^{-2r})^{-\frac{1}{2}}h,
\exp^{-r/2}(1-\frac{R}{2}\exp^{-r} +
\frac{R^2}{12}\exp^{-2r})^{-\frac{1}{2}}Jh\big)
\end{multline*}
and an order $2$ correction term whose effect on the sought formula is zero, see
\cite[proof of lemma 7.6]{ob-mh1} for a precise 
justification of this point. The contribution
of the term involving $R^{op}$ is then
\[ 3\, \left(-\frac{5}{4}\exp^{2r} + \frac{R}{2}\exp^r - \frac{7R^2}{48}\right)
\theta\wedge d\theta ,\]
and that of the last term is
\[ -3\, \left( -\frac{1}{4}\exp^{2r} + \frac{R}{4}\exp^r - \frac{5R^2}{48}\right)
\theta\wedge d\theta .\]
The final boundary term is then
\[ B(\bg,r) = - \frac{d}{8}\exp^{2r} - \frac{\chi}{4}\exp^r + \frac{\chi^2}{24d}\ ,
\]
whose divergent terms cancel exactly those of $3\eta(4\rho^2\theta^2 + \gamma)$,
as expected. Adding the constant terms yield Theorem \ref{maintheo}, and the 
elementary explicit computations that follow it.

\section{Explicitation for lens spaces and proof of the corollaries}

We now specialize the formula obtained in Theorem 1.1 to lens spaces. We will then 
prove the corollaries in a second step.

The lens space $L(p,q)$ is the quotient of the 3-sphere $\cerc^3$ in $\mathbb{C}^2$
by $\ZM/p\ZM$, with its generator acting on $\mathbb{C}^2$ by
$(\exp^{\frac{2i\pi}{p}},\exp^{\frac{2iq\pi}{p}})$, where $q$ is prime
with $p$. They are interesting in connection with filling by Einstein
metrics, since some of them appear as boundary at infinity of seldual
Einstein metrics \cite{cal-sin}. On the other hand, it has been shown
that large families of them admit symplectic fillings. In this case 
our result yields:

\begin{prop}\label{nu:lens}
The $\nu$-invariant of the lens space $L(p,q)$ is
$$\nu(L(p,q))=-\frac{1}{p}+12\, s(p,q,1) . $$
\end{prop}

For sake of comparison, we recall to the interested reader the value of the
classical $\eta$-invariant on lens spaces with the standard round metric, as
computed by Atiyah-Patodi-Singer \cite[Proposition 2.12]{aps2}:
\begin{equation}\label{eq:aps}
\eta(L(p,q)) = -4\, s(p,q,1) .
\end{equation}

\begin{proof}
  For simplicity, we shall assume that $(q-1)$ is prime with $p$ (as a matter
  of fact this implies that we take $q\neq 1$), and we
  leave the general case to the reader.
  Let us see the 3-sphere as the bundle $\mathcal{O}(-1)$ over the
  projective line $\mathbb{C}P^1$. The induced action on 
  $\mathbb{C}P^1$ has two fixed
  points: the two antipodal points, with action of 
  $\ZM/p\ZM$ generated by $\exp^{\pm i2\pi\frac{q-1}{p}}$, and action
  in the fiber by $\exp^{i\frac{2\pi}{p}}$ and
  $\exp^{i2\pi\frac{q}{p}}$ respectively. Therefore
  $L(p,q)$ is a $\cerc^1$-orbifold bundle over an orbifold projective line with two
  orbifold points with angle $\frac{2\pi}{p}$. The Euler
  characteristic is $\chi=\frac{2}{p}$ and the Chern number is
  $d=-\frac{1}{p}$. Now our Theorem \ref{maintheo} and Ouyang's Theorem
  \ref{th:ouy} give the formulas
\begin{align*}
  \nu(L(p,q))&=-3+\frac{2}{p}-12\big(s(p,q-1,1)+s(p,1-q,q)\big) , \\
  \eta(L(p,q))&=1-\frac{1}{p}+4\big(s(p,q-1,1)+s(p,1-q,q)\big) .
\end{align*}
  We deduce $\nu(L(p,q))=-\frac{1}{p}-3\eta(L(p,q))$. The proposition
  then follows from (\ref{eq:aps}).
\end{proof}

For instance, there are large families of lens spaces that admit symplecting 
fillings \cite{lisca}, to which Corollary \ref{cor:3} may be applied. 

\medskip

\noindent\textit{Proof of the corollaries.}
Corollaries \ref{cor:1} and \ref{cor:3}
rely on the formula discovered by the authors \cite[Theorem 1.2]{ob-mh1}:
for any Einstein asymptotically hyperbolic manifold $(M^4,g)$,
\begin{equation} \frac{1}{8\pi^2} \int_M \left( 3|W^-|^2 - |W^+|^2 + \frac{1}{24}
\Scal^2 \right) - \chi(\bar{M}) + 3\,\tau(\bar{M}) = \nu(X).\label{eq:my}
\end{equation}
For complex hyperbolic surfaces, the integral term is zero.
If $\bar M$ is smooth, with $X$ as the only end, then 
the topological contributions always are integers. 
Corollary \ref{cor:1} is then proved.

\medskip

It is instructive to check the results for a holomorphic disk
bundle $D$ over a hyperbolic Riemann surface $\Sigma$, with $X$ as its
boundary. 
Clearly one has $\chi(D)=\chi(\Sigma)=\chi$ and $\tau(D)=-1$.
If $D$ carries a complex hyperbolic metric with $X$ as its boundary at
infinity, then corollary \ref{cor:1} gives the equation
$$ \chi+3 = - \nu(X) = d + 3 + \frac{\chi^2}{4d} $$
and the only solution is $d=\frac{\chi}{2}$.
We then recover the well-known fact that the only disk bundles carrying a
complex hyperbolic metric are the square roots of the (complex) 
tangent bundle.

\medskip

Corollary \ref{cor:3} is again a direct consequence of (\ref{eq:my}),
since for a Kähler-Einstein metric, the integral term is non negative.
For an Einstein metric, the story is more complicated, but positivity
is achieved if solutions to the Seiberg-Witten equations exist, and
it is proven in \cite[corollary 31]{rollin} that it is a consequence of the
nonvanishing of the Kronheimer-Mrowka invariants.

\begin{small}
{\flushleft\sl Acknowledgements}. The authors are grateful to Yoshinobu
Kamishima for useful discussions on possible applications of the
$\nu$-invariant, and to Elisha Falbel for comments.
\end{small}

\bigskip

\bibliographystyle{smfplain}

\medskip

\end{document}